\documentclass[aip,cha,reprint,author-numerical,showkeys]{revtex4-1}
\usepackage{graphicx}
\usepackage{amsmath}
\usepackage{amssymb}
\usepackage{theorem}
\usepackage{dcolumn}
\usepackage{bm}
\usepackage{color}
\usepackage{soul}

\renewcommand{\hl}{}

\numberwithin{equation}{section}

\newtheorem{thm}{Theorem}[section]
{\theorembodyfont{\rmfamily}
\newtheorem{defn}[thm]{Definition}
\newtheorem{conj}[thm]{Conjecture}
\newtheorem{rmk}[thm]{Remark}
}

 \newcommand{\R}{{\mathbb R}}
 \newcommand{\Z}{{\mathbb Z}}

\begin{document}

\title[Conjecture on the Nature of Attractors]{A Test for a Conjecture on the Nature of Attractors for Smooth Dynamical Systems}

\author{Georg A. Gottwald}
\affiliation{School of Mathematics and Statistics, University of Sydney, Sydney 2006 NSW, Australia}
\email{georg.gottwald@sydney.edu.au}
\author{Ian Melbourne}
\affiliation{Mathematics Institute, University of Warwick, Coventry, CV4 7AL, UK}
\email{I.Melbourne@warwick.ac.uk}

\date{\today}

\begin{abstract}
Dynamics arising persistently in smooth dynamical systems ranges from regular dynamics (periodic, quasiperiodic) to strongly chaotic dynamics (Anosov, uniformly hyperbolic, nonuniformly hyperbolic modelled by Young towers). The latter include many classical examples such as Lorenz and H\'enon-like attractors and enjoy strong statistical properties.

It is natural to conjecture (or at least hope) that most dynamical systems fall into these two extreme situations. We describe a numerical test for such a conjecture/hope and apply this to the logistic map where the conjecture holds by a theorem of Lyubich, and to the Lorenz-96 system in $40$ dimensions where there is no rigorous theory.  The numerical outcome is almost identical for both (except for the amount of data required) and provides evidence for the validity of the conjecture.
\end{abstract}


\maketitle
 
\begin{quotation}
A longstanding open problem in the theory of dynamical systems, that continues to be the subject of much discussion by mathematicians and physicists, is the question of what constitutes a typical dynamical system. An answer would not only constitute an immense theoretical advance within the theory of smooth dynamical systems, but would have a profound practical impact on our understanding and analysis of physical phenomena in the real world. In this work we formulate a conjecture on the nature of typical dynamical systems stating that they are either regular or chaotic in a way that assures good statistical properties such as existence of Sinai-Ruelle-Bowen (SRB) measures, exponential decay of correlations, large deviation principles as well as central limit theorems. Since the current state of the theory does not allow for a rigorous theoretical treatment of the conjecture, we devise a numerical test which we use to find corroborating evidence for the conjecture.
\end{quotation}
 
\section{Introduction} 
\label{sec-intro}

A central, but currently intractable, question in the theory of smooth deterministic dynamical systems is to understand the types of attractors for {\em typical} systems.  A classification of attractors would range from very regular dynamics to very chaotic dynamics, including periodic sinks at one extreme and uniformly hyperbolic (Axiom~A) attractors at the other extreme. The uniformly hyperbolic attractors of Smale~\cite{Smale67} generalise the Anosov diffeomorphisms and flows. (Smale's definition of uniformly hyperbolic includes the periodic case, but we shall abuse terminology and reserve the words ``uniformly hyperbolic'' for the nonperiodic case.)

Throughout this paper we are interested in both discrete time dynamical systems (noninvertible maps and diffeomorphisms) and continuous time systems (flows).  Similar comments and results apply to both. However, our notation and definitions will be confined to the discrete case, where $f:\R^n\to\R^n$ is a smooth map with compact attractor $\Lambda\subset\R^n$. Our focus is primarily on dissipative systems, but the material goes over to Hamiltonian systems with the obvious modifications.

An important property of uniformly hyperbolic attractors is the existence of a {\em physical measure}, or {\em SRB measure} after Sinai, Ruelle and Bowen, which has the property that time averages converge to the space average for a set of initial conditions of positive volume (ie.\ positive Lebesgue measure).  This is in contrast to the ergodic theorem for ordinary ergodic measures where the convergence takes place for a set that has full measure with respect to the ergodic measure which however is usually a set of zero volume (since the ergodic measure is supported on the attractor $\Lambda$ which is usually of zero volume).

\begin{defn} An ergodic measure $\mu$ supported on $\Lambda$ is an {\em SRB measure} if  there is a set $B$ of positive volume such that
\[
\lim_{n\to\infty}\frac1n\sum_{j=0}^{n-1}v(f^jx_0)=\int_\Lambda v\,d\mu,
\]
for every continuous observable $v:\R^n\to\R$ and for all $x_0\in B$.
\end{defn}

Uniformly hyperbolic attractors have numerous strong statistical properties. In particular, they have exponential decay of correlations up to a finite cycle~\cite{Sinai72,Bowen75,Ruelle78}.

\begin{defn}
An attractor $\Lambda$ with ergodic measure $\mu$ has {\em exponential decay of correlations} if there exists a constant $\gamma\in(0,1)$ such that for all smooth $v,w:\R^n\to\R$ there is a $C>0$ such that
\[
\Bigl|\int_\Lambda v\,w\circ f^n\,d\mu
-\int_\Lambda v\,d\mu
\int_\Lambda w\,d\mu\Bigr|\le C\gamma^n.
\]
More generally, the attractor has {\em exponential decay of correlations (up to a finite cycle)} if there exists $k\ge1$ and disjoint compact sets $\Lambda_1,\dots,\Lambda_k$ such that for all $i=1,\ldots, k$ it is the case that $f(\Lambda_i)=\Lambda_{i+1}$ (with $k+1=1$) and $f^k:\Lambda_i\to\Lambda_i$ has exponential decay of correlations.
\end{defn}

From now on, we omit the words ``up to a finite cycle'' and speak simply of exponential decay of correlations.   The SRB measure $\mu$ on a uniformly hyperbolic attractor enjoys this property and it suffices that $v$ and $w$ are Lipschitz (or even H\"older, in which case the constant $\gamma$ depends on the H\"older class). There are numerous other statistical properties such as central limit theorems that hold for uniformly hyperbolic attractors.  These are described in a more general setting below.

\begin{rmk} Decay of correlations for uniformly hyperbolic (even Anosov) flows is rather less well understood.  Only partial results exist~\cite{Dolgopyat98a,Dolgopyat98b,Liverani04,FMT07}; however statistical limit laws such as central limit theorems and invariance principles remain valid for uniformly hyperbolic flows~\cite{Ratner73,DenkerPhilipp84,MN09}.
\end{rmk}

Smale conjectured (see for example Section 2.5 in Palis~\cite{Palis05}) that for typical dynamical systems (typical in the sense of $C^r$ open and dense, $r\ge1$) periodic sinks and uniformly hyperbolic attractors comprise the full range of possibilities.  This conjecture turned out to be false, and moreover the notion of typicality turned out to be inadequate even in situations where the conjecture holds (see for example items (i)--(iii) below). Over the last 40--50 years, numerous examples have arisen that make it necessary to enlarge the notions of being very regular or very chaotic.

\begin{itemize}
\item[(i)] KAM tori with quasiperiodic dynamics  are nonrobust in a topological sense (they are destroyed by $C^r$ small perturbations,) but they are unavoidable in a probabilistic sense (the set of parameters that give rise to KAM tori has large measure).   For dissipative systems a similar phenomenon arises in Naimark-Sacker bifurcation from a periodic solution.
\item[(ii)]  The logistic map (see Section~\ref{sec-log} for more details) is a one-parameter family of one-dimensional maps. For each value of the parameter there is a unique attractor that attracts almost every trajectory.  For an open and dense set of parameters, the attractor is a periodic sink. However, Jakobson~\cite{Jakobson81} showed that the complementary set of parameters has positive measure. More recently, Lyubich~\cite{Lyubich02} proved that almost every parameter in this complementary set satisfies the so-called Collet-Eckmann condition~\cite{ColletEckmann83} and hence constitutes strongly chaotic (though not uniformly hyperbolic) dynamics.
\item[(iii)]  H\'enon-like attractors~\cite{Henon76} arise near quadratic homoclinic tangencies~\cite{BenedicksCarleson91,MoraViana93} and are strongly chaotic~\cite{BenedicksYoung00}.  These are again unavoidable in a probabilistic sense.
\item[(iv)]  Geometric Lorenz attractors are topologically robust but nonuniformly hyperbolic examples of strongly chaotic systems~\cite{AfraimovicBykovSilnikov77,GuckenWilliams79,Williams79}.  Tucker~\cite{Tucker02} showed that these include the classical Lorenz attractor~\cite{Lorenz63}.
\end{itemize}

The strongly chaotic attractors mentioned above -- uniformly hyperbolic, Collet-Eckmann, H\'enon-like, (geometric) Lorenz -- have the common property that they are modelled by a Young tower with exponential tails as introduced by Young~\cite{Young98}.  (For Lorenz attractors, it is the Poincar\'e map that is modelled by a Young tower.) Roughly stated, a dynamical system $f:\Lambda\to\Lambda$ is modelled by a Young tower with exponential tails if there exists a set $Y\subset\Lambda$ with return time function $\tau:Y\to\Z^+$ (not necessarily the first return time) and return map $F=f^\tau : Y \to Y$ such that (i) $F$ is uniformly hyperbolic, and (ii) the likelihood of a large return time $\tau$ is exponentially small.

Numerous strong statistical properties have been proved for such attractors modelled by Young towers: existence of an SRB measure, exponential decay of correlations and central limit theorems~\cite{Young98}, large deviation principles~\cite{MN08,ReyBelletYoung08}, Berry-Ess\'een estimates and local limit theorems~\cite{Gouezel05}, invariance principles~\cite{MN05,MN09,Gouezel10}. There is also an enlarged class of attractors~\cite{Young99} that possess polynomial decay of correlations; where this decay is summable the above statistical properties apply.

In a sense that can be made precise, there is an equivalence between the existence of a Young tower and strong statistical properties~\cite{AFLV11}. This observation uses the work of Alves {\emph{et al.}}~\cite{AlvesLuzzattoPinheiro04} and Melbourne \& Nicol~\cite{MN08}.

Since there are good reasons for hoping (if not believing) that most attractors are either highly regular or enjoy strong statistical properties, and in the absence of convincing counterexamples, one possibility is to define strongly regular attractors to be the periodic and quasiperiodic ones, and strongly chaotic attractors to be the ones modelled by a Young tower with exponential tails.  This leads naturally to the following deliberately imprecise conjecture.

\begin{conj} \label{conj}
Typically (in a sense that we do not make precise), the attractors for smooth dynamical systems fall into one of the following two classes:
\begin{itemize}
\item[(a)]
Regular dynamics: $\Lambda$ is a periodic or quasiperiodic sink.
\item[(b)]
Chaotic dynamics: $\Lambda$ is modelled by a Young tower with exponential tails.
\end{itemize} 
(In the case of flows, this statement is at the level of the Poincar\'e map.)
\end{conj}

A precise conjecture would require a precise definition of ``typically'', probably leading to the failure, though not necessarily the relevance, of the conjecture.

\subsubsection*{A test for Conjecture~\ref{conj}}

Although it is hard to see how to test directly for Conjecture~\ref{conj}, there are certain implications that can be tested numerically.  Suppose that $\Lambda$ is an attractor for a map or diffeomorphism $f:\R^n\to\R^n$ and that $\mu$ is an ergodic invariant measure on $\Lambda$.  Let $v:\R^n\to\R$ be a smooth observable. Recall that the power spectrum $S:[0,2\pi]\to[0,\infty)$ is given by
\[
S_\omega=\lim_{n\to\infty}\frac1n S_\omega(n), \quad
S_\omega(n)=\int_\Lambda\Bigl|\sum_{j=0}^{n-1} e^{ij\omega} v\circ f^j\Bigr|^2\,d\mu.
\]
Since $S_{2\pi-\omega}=S_\omega$ we restrict from now on to the interval $[0,\pi]$.

The following dichotomy was established by Melbourne \& Gottwald~\cite{MG08}.
\begin{thm}   \label{thm-power}
Let $\Lambda$ be a periodic or quasiperiodic sink, or an attractor modelled by a Young tower with exponential tails for a smooth map $f:\R^n\to\R^n$.  Suppose that $\mu$ is the SRB measure on $\Lambda$.  Let $v:\R^n\to\R$ be a $C^\infty$ observable. Typically, 
\begin{itemize}
\item[(a)]  In the periodic/quasiperiodic case, $S_\omega=0$ almost everywhere.
Moreover,
\[
K_\omega=\lim_{n\to\infty} \log S_\omega(n)/\log n=0
\]
for all but finitely many $\omega\in[0,\pi]$.
\item[(b)]  In the Young tower case, there is a constant $s_0>0$ such that $S_\omega\ge s_0$ for all but finitely many values of $\omega$. In particular, 
\[
K_\omega=\lim_{n\to\infty} \log S_\omega(n)/\log n=1
\]
for all but finitely many $\omega\in[0,\pi]$.
\end{itemize}
\end{thm}

\begin{rmk}
Often in the physics literature, regular and chaotic dynamics is distinguished in terms of the power spectrum~\cite{GollubSwinney75}.  Broadband power spectrum (where there exists an interval or at least a set of positive measure on which $S_\omega$ is positive) is seen as being the signature of chaotic dynamics.

The dichotomy in Theorem~\ref{thm-power} is significantly stronger.   In case~(a), we are requiring that $S_\omega(n)$ grows slower than any polynomial rate.  In contrast, the requirement that $S_\omega=0$ is compatible with values of $K_\omega$ anywhere in $[0,1]$. For example, if $S_\omega(n)$ grows like $n/\log n$, then $S_\omega=0$ but $K_\omega=1$.

In case~(b), the power spectrum is positive and bounded away from zero for all but finitely many points, which is rather more than claiming broadband spectrum. 
\end{rmk}

Based on Theorem~\ref{thm-power}, our conjecture can be tested as follows. Consider a parameterized family of smooth dynamical systems, with parameter $a\in\R$.   For each fixed value of $a$, compute the family of limits $S_\omega$ and check whether they are almost all zero or almost all one.  Such a test can be carried out numerically by taking values of $a$ and $\omega$ that are reasonably dense and estimating the growth rate of $S_\omega(n)$.

\begin{rmk} 
A slightly weaker version of the conjecture would be to include Young towers with polynomial tails (rather than only those with exponential tails).    Our test does not distinguish between these situations.   However, we do not know of any persistent examples in smooth dynamics where an attractor is modelled by a Young tower with subexponential tails, but not by a Young tower with exponential tails.
\end{rmk}

There are similarities and differences between the test proposed above and the $0$--$1$ test for chaos~\cite{GM04,GM05,GM09,GM09b}.  The $0$--$1$ test is optimised to work with limited amounts of data.   In particular, taking the median value of $K_\omega$ for $100$ randomly chosen values of $\omega$ greatly accelerates the convergence of the test.   The test described in this paper is a much more stringent examination of the dichotomy in Conjecture~\ref{conj} but requires much more data.  Even for the logistic map,
the refined test in this paper requires enormous amounts of data that would be impractical in the $0$--$1$ test for chaos.   

\vspace{1ex}

Our conjecture is related to the Palis conjecture and to the Gallavotti-Cohen chaotic hypothesis.

\vspace{.3ex}
\noindent{\bf Palis conjecture.}  
As already mentioned, Smale's conjecture regarding the ubiquity of periodic sinks and uniformly hyperbolic attractors turned out to be false.  Hence it became necessary to formulate a weakened statement.  Over the years, Palis gave a number of conjectures in this direction; we refer to the original work by Palis~\cite{Palis00,Palis05,Palis08} for statements of these conjectures (rather more precise than ours!) and progress towards their verification.  

The emphasis in the Palis conjectures is global, focusing on (i) the finitude of attractors possessing SRB measures, with the property that the union of their basins accounts for a set of initial conditions of full measure, and (ii) the stability of these attractors under perturbations. Our conjecture is more local since we have said nothing about the finitude of attractors, nor their stability under perturbations.  However, for typical attractors taken on their own, we make stronger statements about their statistical properties.

\vspace{.3ex}
\noindent{\bf Gallavotti-Cohen chaotic hypothesis.}
The chaotic hypothesis~\cite{GallavottiCohen95} proposes that chaotic systems should be considered as Anosov systems for practical purposes. Since the property of part (b) of Theorem~\ref{thm-power} is certainly valid for Anosov systems, our numerical test can be viewed as a test also of the chaotic hypothesis.

\vspace{1ex}
The remainder of the paper is organised as follows. In Section~\ref{sec-test}, we describe how to test parametrized families of dynamical systems for their conformity to Conjecture~\ref{conj}. In Section~\ref{sec-log}, we carry out this test for the logistic map. This provides a benchmark for our test since the conjecture is known to be valid by Lyubich~\cite{Lyubich02}. In Section~\ref{sec-96}, we carry out the test for the $40$-dimensional Lorenz-96 system, which is regarded as highly important in meteorological studies, and which is far beyond the current understanding of rigorous dynamical systems theory.  Nevertheless, the numerical results for Lorenz-96 are similar to those for the logistic map except for the amount of data required for convergence. We conclude with a brief summary in Section~\ref{sec-discussion}.

\section{The numerical test}
\label{sec-test}

Consider a smooth family of maps $f_a:\R^n\to\R^n$ where $a\in\R$ is a parameter.  For convenience, we assume that all trajectories are bounded.

Suppose that $a\in\R$ and that $\Lambda\subset\R^n$ is an attractor for $f_a$. For values of $\omega$ chosen randomly from $[0,\pi]$ we compute $K_\omega$ as defined in Theorem~\ref{thm-power}. Then according to Conjecture~\ref{conj} and Theorem~\ref{thm-power} we anticipate that $K_\omega$ takes the constant value $0$ or $1$ independent of $\omega$ (for all but finitely many $\omega$).

To carry out this procedure numerically, we note that computing $S_\omega(n)$ directly is unfeasible since $\Lambda$ (and $\mu$) are not given. However, by the ergodic theorem, for $\mu$-almost every $x_0\in\Lambda$,
\begin{align} \label{eq-MSQ}
S_\omega(n)=\lim_{J\to\infty}\frac1J\sum_{j=0}^{J-1}|p_\omega(j+n)-p_\omega(j)|^2,
\end{align}
where 
\begin{align} \label{eq-p}
p_\omega(n)=\sum_{\ell=0}^{n-1}e^{i\ell\omega}v(f_a^\ell x_0).
\end{align}
Moreover, assuming the conjecture, typically $\mu$ can be taken to be an SRB measure and $x_0$ can be chosen from a set of positive Lebesgue measure.

Equally, assuming the conjecture, typically $x_0\in\R^n$ lies in the basin of an attractor $\Lambda$ (depending on $x_0$ and $a$) and $S_\omega(n)$ can be computed using~\eqref{eq-MSQ} and~\eqref{eq-p}. From this $K_\omega=\lim_{n\to\infty}\log S_\omega(n)/\log n$ can be computed. Again, if Conjecture~\ref{conj} is valid, then for typical $f_a$ and $x_0$ it should be the case that $K_\omega$ takes the constant value either $0$ or $1$ independent of $\omega$ (for all but finitely many $\omega$).

There are some finite computation issues that need to be addressed. The most crucial one is that the definition of $K_\omega$ involves a double limit, first as $J\to\infty$ and then as $n\to\infty$.  We will ignore this issue to begin with and return to it at the end of the section.

To implement the test, we take as initial condition $x_0=f^{1000}x_1$ where $x_1$ is chosen at random and fixed throughout.  (Neglecting this transient of $1000$ iterates is not strictly necessary but speeds up the calculations.) A finite but reasonably dense set of parameters $a$ is specified. For each value of $a$, we compute $K_\omega$ for $100$ (say) randomly chosen values of $\omega\in[0,\pi]$.  Given the finiteness of the data, it is necessary to specify small open intervals $I_0$ and $I_1$ containing $0$ and $1$ respectively, such that $K_\omega\in I_r$ is viewed as an $r$ for $r=0,1$.  
Define
\begin{align*}
& M_0=\#\{\omega:K_\omega\in I_0\},\\
& M_1=\#\{\omega:K_\omega\in I_1\},\\
& M_u=\#\{\omega:K_\omega\not\in I_0\cup I_1\},
\end{align*}
with $M_0+M_1+M_u=100$. (Here $u$ stands for undecided.) As $K_\omega$ is computed with greater and greater precision, a consequence of the conjecture is that either $M_0\to100$ or $M_1\to100$. An implication that is easier to test for is that 
\[
M_u\to0  \quad\text{and} \quad \min\{M_0,M_1\}\to0.
\]

The numerical test that we propose can now be stated more precisely. We make three choices of intervals  
\begin{itemize}
\item[(i)] $I_0=(-0.1,0.3)$, $I_1=(0.7,1.1)$.
\item[(ii)] $I_0=(-0.1,0.2)$, $I_1=(0.8,1.1)$.
\item[(iii)] $I_0=(-0.1,0.1)$, $I_1=(0.9,1.1)$.
\end{itemize}
For each of these choices, we take $A$ equally spaced values of the parameter $a$ and $100$ values of $\omega\in[0,\pi]$ chosen at random. (The value of $A$ will depend on the length of the range of interesting parameters for the dynamical system.) Then we analyse the convergence to zero of the following four quantities as the limit $J\to\infty$ and $n\to\infty$ is approached:
\begin{align*}
& Q_u  ={\textstyle\sum_a} M_u,\\
& Q_u'=\#\{a:M_u>10\},
\end{align*}
and
\begin{align*}
& Q_{\rm min}={\textstyle\sum_a}\min\{M_0,M_1\}, \\
& Q_{\rm min}'=\#\{a:\min\{M_0,M_1\}>10\}.
\end{align*}
The quantity $Q_u\in\{0,1,\dots,100A\}$ denotes the total number of values of $\omega$ and parameter values $a$ for which the value of $K$ is undecided (i.e.\ it lies outside $I_0\cup I_1$), whereas the quantity $Q_u'\in\{0,1,\dots,A\}$ denotes the number of parameter values for which $K_\omega$ is undecided for more than $10\%$ of the choices of $\omega$. Similarly for $Q_{\rm min}$ and $Q'_{\rm min}$ with the number of undecideds replaced for each $a$ by the minimum of the number of $0$s and the number of $1$s.

\hl{Our choices for the intervals $I_0$ and $I_1$ around $0$ and~$1$ are somewhat arbitrary; if the conjecture is true than eventually the four quantities $Q_u$, $Q_u'$, $Q_{\rm min}$,
$Q_{\rm min}'$ will reach zero, regardless.  The amount of data to achieve this depends on the choice of intervals, but this dependence is not relevant for Conjecture}~\ref{conj}.  \hl{On the other hand, our numerical experiments indicate that $K_\omega$ very quickly lies between $0$ and $1$ (within a small error), and we have chosen the sharper lower limit $-0.1$ for $I_0$ and
upper limit $1.1$ for $I_1$ with this in mind.}

\hl{We note that the convergence to zero need not be monotone.  
For example, suppose that a parameter value $a$ yields a chaotic attractor of
class (b), so that eventually $M_0=0$ and $M_1=100$.  If the convergence is sufficiently slow, then it is possible that $M_0=95$ and $M_1=5$ (say) for $N$ too small.  For moderate values of $N$, the situation might improve to $M_0=15$, $M_1=85$.   In this case, the parameter $a$ contributes adversely to $Q_{\rm min'}$ for $N$ moderate but not for $N$ small.
An example of this is shown in Figure}~\ref{fig-logistic} \hl{where the number of outliers for $Q_{\rm min}'$ increases from zero to one as $N$ increases within the range of our experiment.}

\subsubsection*{The double limit}

As promised, we discuss the issues regarding the double limit in the formula 
\[
K_\omega=\lim_{n\to\infty}\lim_{J\to\infty}\log\Bigl(\frac1J 
\sum_{j=0}^{J-1}|p_\omega(j+n)-p_\omega(j)|^2\Bigr)
/\log n.
\]
Under certain conditions, it should be possible to prove that for any $\tau\in(0,1)$,
\begin{align} \label{eq-g}
K_\omega=\lim_{J\to\infty}\log\Bigl(\frac1J
\sum_{j=0}^{J-1}|p_\omega(j+J^{\,\tau})-p_\omega(j)|^2\Bigr)
/\log J^{\,\tau}.
\end{align} 
However, there is no way to tell how large $J$ needs to be for a given $\tau$ to be effective, rendering formula~\eqref{eq-g} unsuitable for a numerical test.  We follow the simpler route of replacing $J^{\tau}$ by  $\delta J$ where $\delta$ is a small constant (depending on the family of dynamical systems).  By inspection for a few randomly chosen values of $a$ and $\omega$, we check that $\delta$ is sufficiently small for the range of $n$ used in the numerical test.   In Sections~\ref{sec-log} and~\ref{sec-96}, we verify that $\delta=0.01$ suffices 
for the logistic map and the Lorenz-96 system, respectively.

Suppose that $N$ denotes the number of iterates available for the numerics, so we have computed $f_a^jx_0$ for $j=0,\dots,N-1$. Then $S_\omega(n)$ can be computed for $n=0,1,\dots,\delta N$ and we can use $\log S_{\omega}(\delta N)/\log \delta N$ as an estimate for $K_\omega$.

\subsubsection*{Speeding up the test}

We have already mentioned that taking a short transient (say $1000$ iterates) speeds up the convergence in the test.    There are further devices for speeding up the test that we observed while developing the $0$--$1$ test for chaos.

First, it is useful to define the modified mean square displacement~\cite{GM09,GM09b}
\[
D^0_\omega(n)=S_\omega(n)-({\textstyle\int}_\Lambda v\,d\mu)^2\frac{1-\cos n\omega}{1-\cos\omega}.
\]
Note that $\int_\Lambda v\,d\mu=\lim_{n\to\infty}\frac1n\sum_{j=0}^{n-1}v(f^jx_0)$ can be computed using the ergodic theorem, and that $K_\omega=\lim_{n\to\infty}\log D^0_\omega(n)/\log n$. For Young towers with exponential tails that are mixing, it was proved that the convergence as $n\to\infty$ is now uniform in $\omega$~\cite{GM09b}.  Even in the nonmixing case, numerics~\cite{GM09b} show this to be a useful modification.  To avoid taking logarithms of negative numbers, we set $D_\omega(n)=D^0_\omega(n)+C$ where $C=\max_{k=1,\dots,\delta N}|D^0_\omega(k)|$. \hl{(Since $C$ is a constant once the amount of data is specified, the growth rate in $n$ is unchanged.)}  Then we replace $S_\omega(n)$ by $D_\omega(n)$ in all the formulas.

Second, to make more efficient use of the data, we compute $D_\omega(n)$ for all $n\le \delta N$ and perform linear regression on $\log D_\omega(n)$ plotted against $\log n$ (cf.\ Gottwald \& Melbourne~\cite{GM04,GM09}).

\section{Logistic map}
\label{sec-log}

The logistic map, or quadratic family, $f_a(x)=a x(1-x)$, $0\le a\le4$, is a convenient example to begin with since it is very well understood. For each value of $a$, there is a unique attractor $\Lambda\subset[0,1]$. For an open and dense set of parameters $a\in P\subset[0,4]$, the attractor $\Lambda$ is a periodic sink. However, Jakobson~\cite{Jakobson81} (see also Benedicks \& Carleson~\cite{BenedicksCarleson85}) proved that ${\rm Leb}(P)<4$. By Lyubich~\cite{Lyubich02} for almost every $a\in [0,4]\setminus P$, the Collet-Eckmann condition~\cite{ColletEckmann83} holds, and this implies the existence of a Young tower with exponential tails (see Theorem 7 in Young~\cite{Young98}). Hence Conjecture~\ref{conj} is valid for this family.

We confine our numerics to the parameter range $3.5\le a\le 4$ since $[0,3.5]\subset P$ and corresponds entirely to periodic sinks of low period (at most period $4$).

\begin{figure}[h]
\includegraphics[width=0.8\columnwidth]{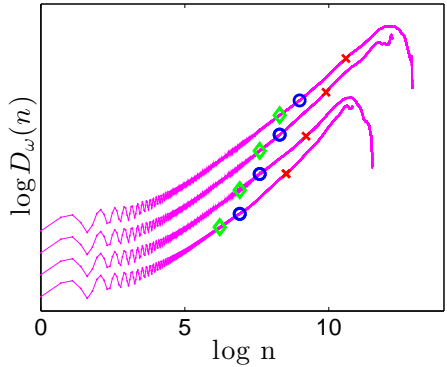}
\caption{Graph of $\log D_\omega(n)$ against $\log n$ for the logistic map with $a=3.6022$ and $\omega=1.9418$.   We take $n=1,\dots,N$ with $N=50,000$, $N=100,000$, $N=200,000$, $N=400,000$.  The various ranges used for the linear regression are marked by green diamonds for $\delta=0.01$, blue circles for $\delta=0.02$ and red crosses for $\delta=0.1$.
\hl{The four graphs are spaced apart vertically so that they can be seen separately, with $N=50,000$ at the bottom, up to $N=400,000$ at the top.}}
\label{fig-logistic_delta}
\end{figure}

The first step is to determine a suitable value of $\delta$. To achieve this, we chose various values of $a\in[3.5,4]$ and $\omega\in[0,2\pi]$ at random and plotted $\log D_\omega(n)$ against $\log n$ for various ranges of \hl{$n=1,\dots,N$}. Theorem~\ref{thm-power} implies that the graph should be linear, but in practice given $N$ iterates of the dynamical system, the graph is linear only up to a certain point.   A typical example is shown in Figure~\ref{fig-logistic_delta}.  \hl{(The graphs for different choices of $N$ need not coincide for a given $n$ since the averaging in}~\eqref{eq-MSQ} \hl{is over a different range.
The strange (and inconsistent) behaviour for large $n$ confirms that there is insufficient averaging once $n$ is too large relative to $N$.)}
It is evident from the graphs that $\delta=0.1$ is too large, whereas $\delta=0.01$ is comfortably within the linear range. Our experiments with various choices of $a$ and $\omega$ confirm that $\delta=0.01$ is a safe choice for the entire range of values of $a$, $\omega$ and $N$ in our numerical test.  From now on we fix this value of $\delta$ for the logistic map.

\begin{figure*}[htb]
\includegraphics[width=0.8\columnwidth]{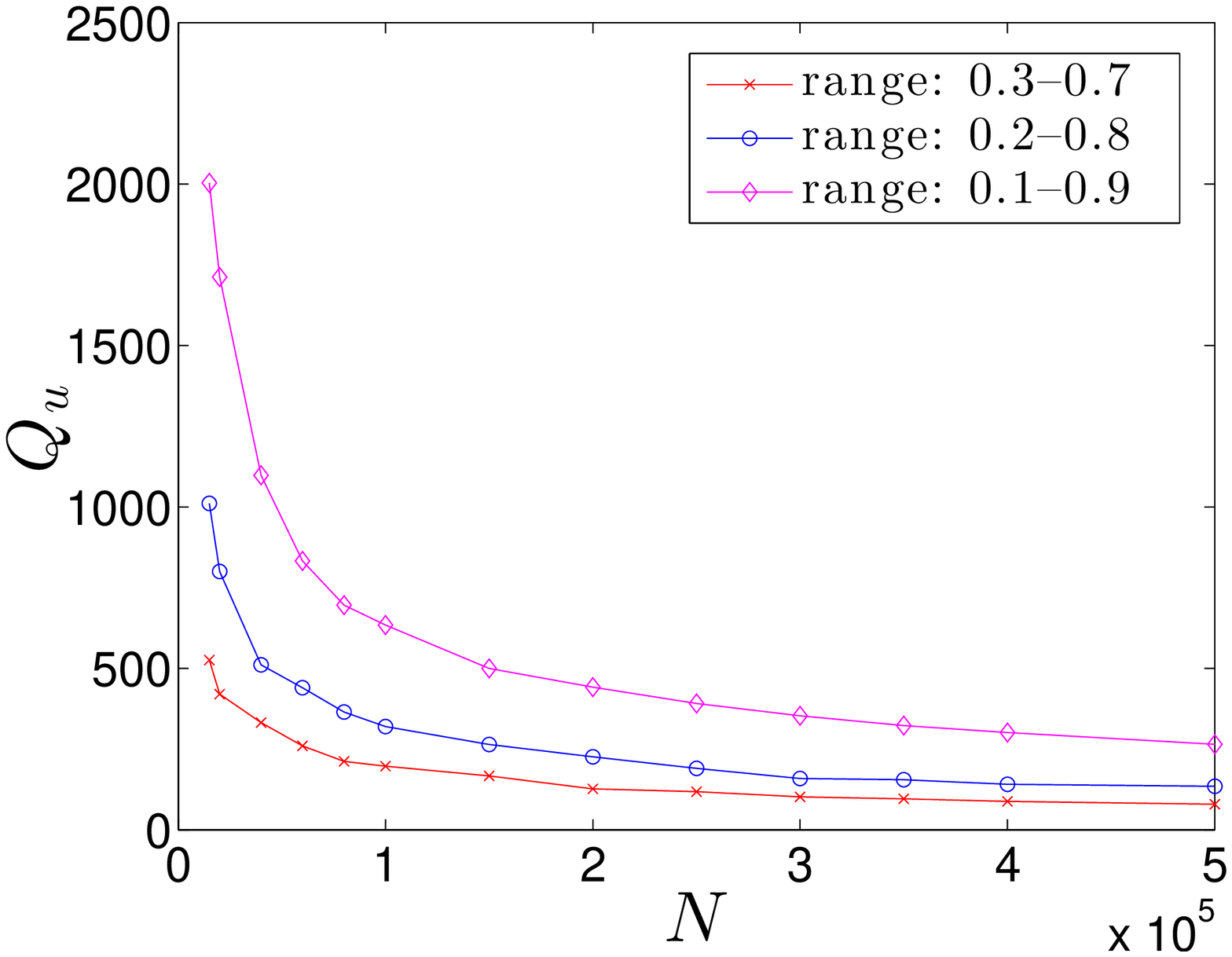}
\includegraphics[width=0.8\columnwidth]{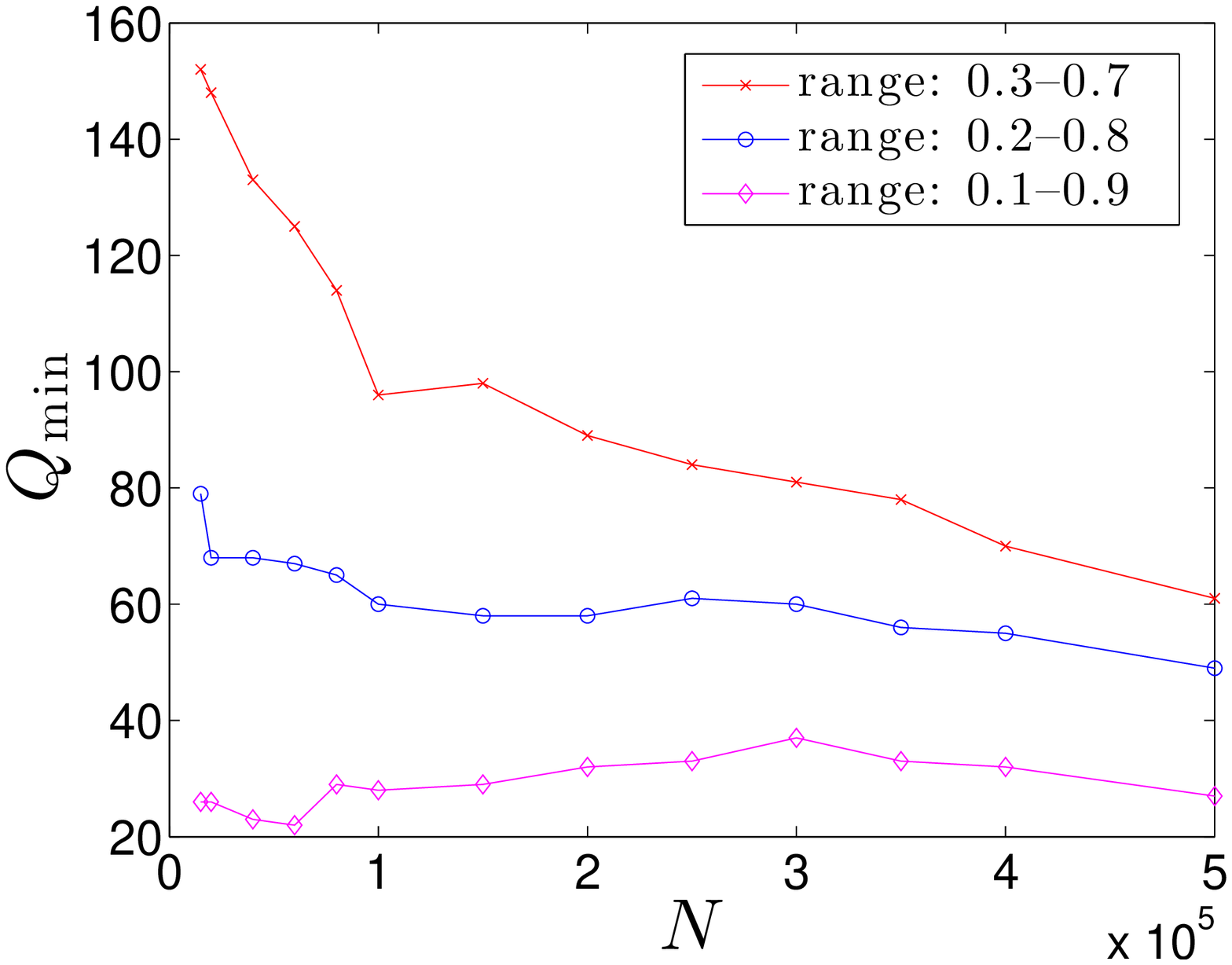}
\\
\includegraphics[width=0.8\columnwidth]{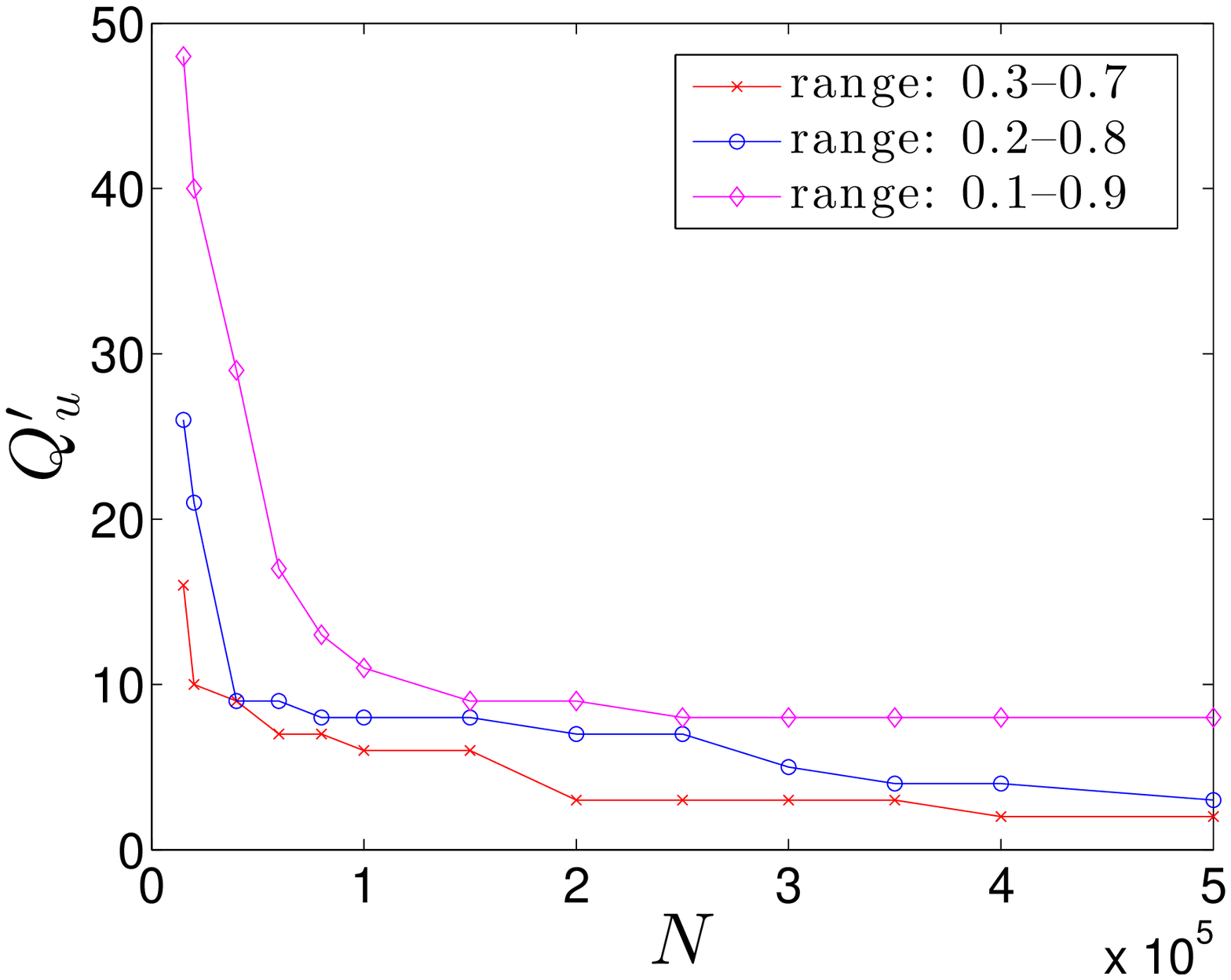}
\includegraphics[width=0.8\columnwidth]{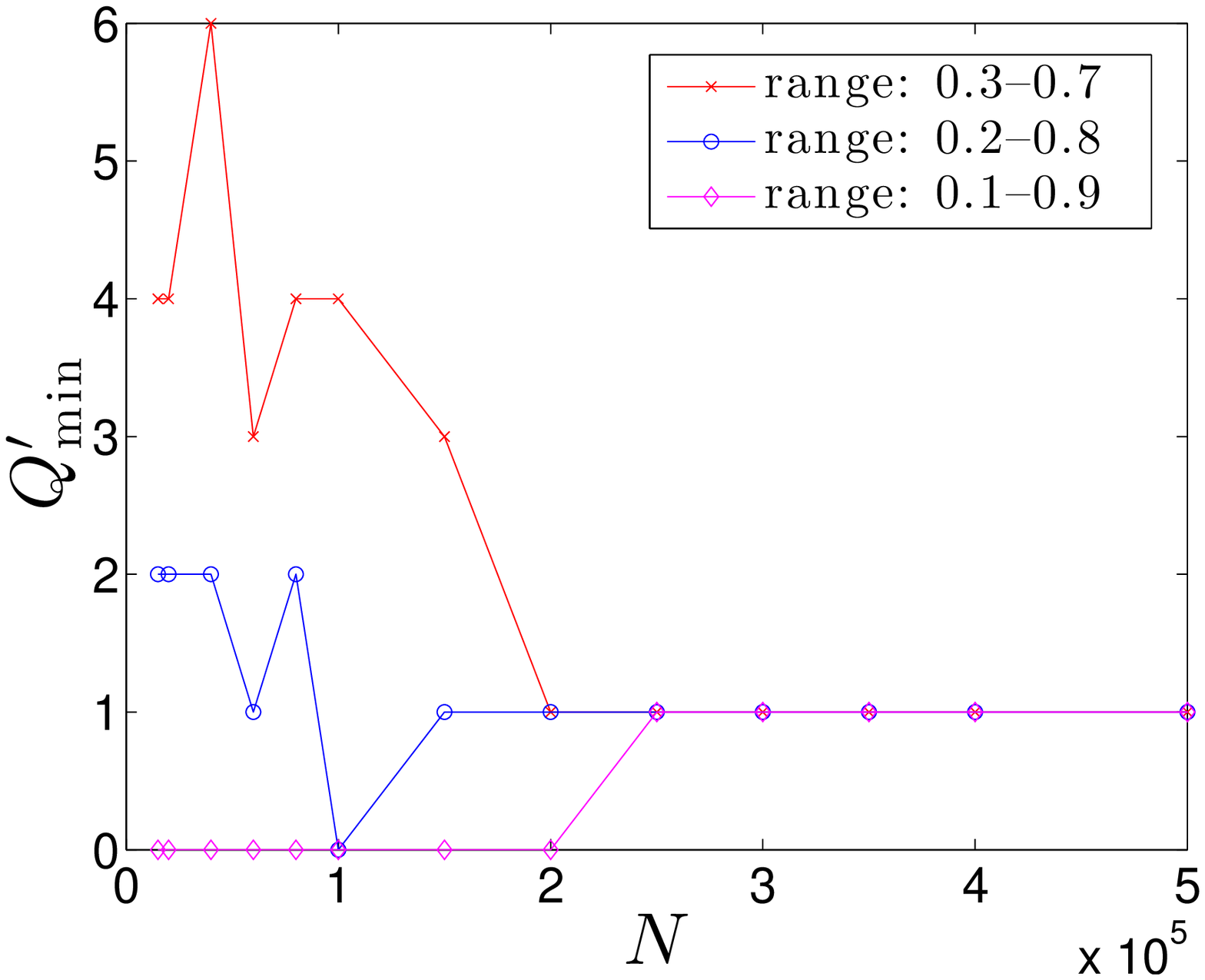}
\caption{Graphs of $Q_u$, $Q_{\rm min}$, 
$Q_u'$, $Q_{\rm min}'$ against the number of iterates $N$ for the logistic map. In each case the results are shown for the three choices 
(i) $I_0=(-0.1,0.3)$, $I_1=(0.7,1.1)$, red crosses;
(ii) $I_0=(-0.1,0.2)$, $I_1=(0.8,1.1)$, blue circles;
(iii) $I_0=(-0.1,0.1)$, $I_1=(0.9,1.1)$, magenta diamonds.
}
\label{fig-logistic}
\end{figure*}

Our numerical results for the logistic map are shown in Figure~\ref{fig-logistic}. The results are consistent with the theory, based on Lyubich~\cite{Lyubich02}, which dictates that the four quantities $Q_u$, $Q_{\rm min}$, $Q'_u$, $Q'_{\rm min}$ converge to zero as $N\to\infty$.  However, it is also clear that there are a handful of cases that are converging very slowly, with little appreciable improvement from $N=100,000$ to $N=500,000$. 
\hl{It is well-known that the onset of chaos near $a\approx 3.57$ leads to very slow convergence in any numerical method for distinguishing regular and chaotic dynamics. Nevertheless, by}~\cite{Lyubich02} \hl{we know that the conjecture is true for this example, so the difficulty is not with the conjecture itself, but with the numerical verification of the conjecture.}
Understanding these limitations to \hl{this (or any)} numerical test is instructive when applying it to examples where there is no proof of convergence.

\begin{figure}[htb]
\includegraphics[width=0.8\columnwidth]{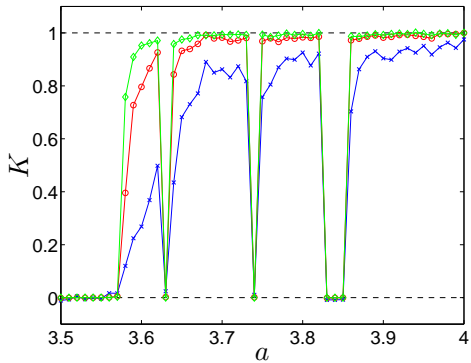}
\caption{The median value of $K(\omega)$ plotted against the parameter
$a$ for the logistic map.  
Blue crosses: $N=10,000$ iterates. 
Red circles: $N=100,000$ iterates. 
Green diamonds: $N=500,000$ iterates.}
\label{fig-logistic_01}
\end{figure}

To this end, it is useful to first contrast these results with the $0$--$1$ test for chaos which uses the median value of $K(\omega)$, and hence converges very quickly for most values of the parameter $a$, see Figure~\ref{fig-logistic_01}.   The problematic parameters are \hl{indeed} the ones near the onset of chaos $a\approx 3.57$ and \hl{also} near the first periodic window $a\approx 3.63$.  It is noteworthy that the onset of chaos after the large period $3$ window near $a\approx 3.83$ does not cause a problem.  (Of course, periodic windows are dense but at this level of resolution, where $a$ is increased in increments of $0.01$, there are only three periodic windows.)

\begin{figure}[htb]
\includegraphics[width=0.8\columnwidth]{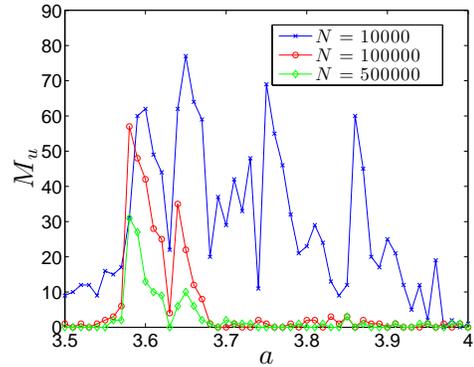}
\caption{Graph of the percentage $M_u$ of undecided values of $\omega$ plotted against the parameter $a$ 
for the logistic map, using the range
(ii) $I_0=(-0.1,0.2)$, $I_1=(0.8,1.1)$ throughout.
Blue crosses: $N=10,000$ iterates. 
Red circles: $N=100,000$ iterates. 
Green diamonds: $N=500,000$ iterates.}
\label{fig-logistic_outlier}
\end{figure}

It is easily verified that the eventually slow convergence in Figure~\ref{fig-logistic} is entirely connected with the rare problematic parameters indicated in Figure~\ref{fig-logistic_01}. In particular, a relatively small number of outliers persist for a very large number of iterates.  This is illustrated in Figure~\ref{fig-logistic_outlier}, where the number of undecided values $M_u$ is shown as a function of the parameter $a$ for a range of iterates $N$. It is shown that by the time we reach $N=500,000$ iterates, the nonconvergence of the four quantities in Figure~\ref{fig-logistic} is due almost entirely to two values of the parameter, $a=3.58$ and $a=3.59$.  

\section{Lorenz-96 model}
\label{sec-96}

In this section, we consider the Lorenz-96 model
\begin{align}
\label{lorenz96}
\frac{dx_i}{dt}=x_{i-1}(x_{i+1}-x_{i-2})-x_i+a \quad {\rm{with}} \quad
i=1,\cdots,m\; ,
\end{align}
where $x_0=x_m$. This system of ordinary differential equations was first introduced by Lorenz as an idealised model for midlatitude atmospheric dynamics \cite{Lorenz96,LorenzEmanuel98}. We consider the case $m=40$ with the parameter $a$ varying in the interval $[3,7]$ in increments of $0.1$. Throughout, we integrate the system using a time step of $0.0005$ and record, after an initial transient of $10,000$ time steps, $N$ data points after each $1000$ time steps. As an observable we take $\phi=x_1$, so $\phi(n)=x_1(t)$ with $t=0.5\, n$.\\

\begin{figure}[h]
\includegraphics[width=0.8\columnwidth]{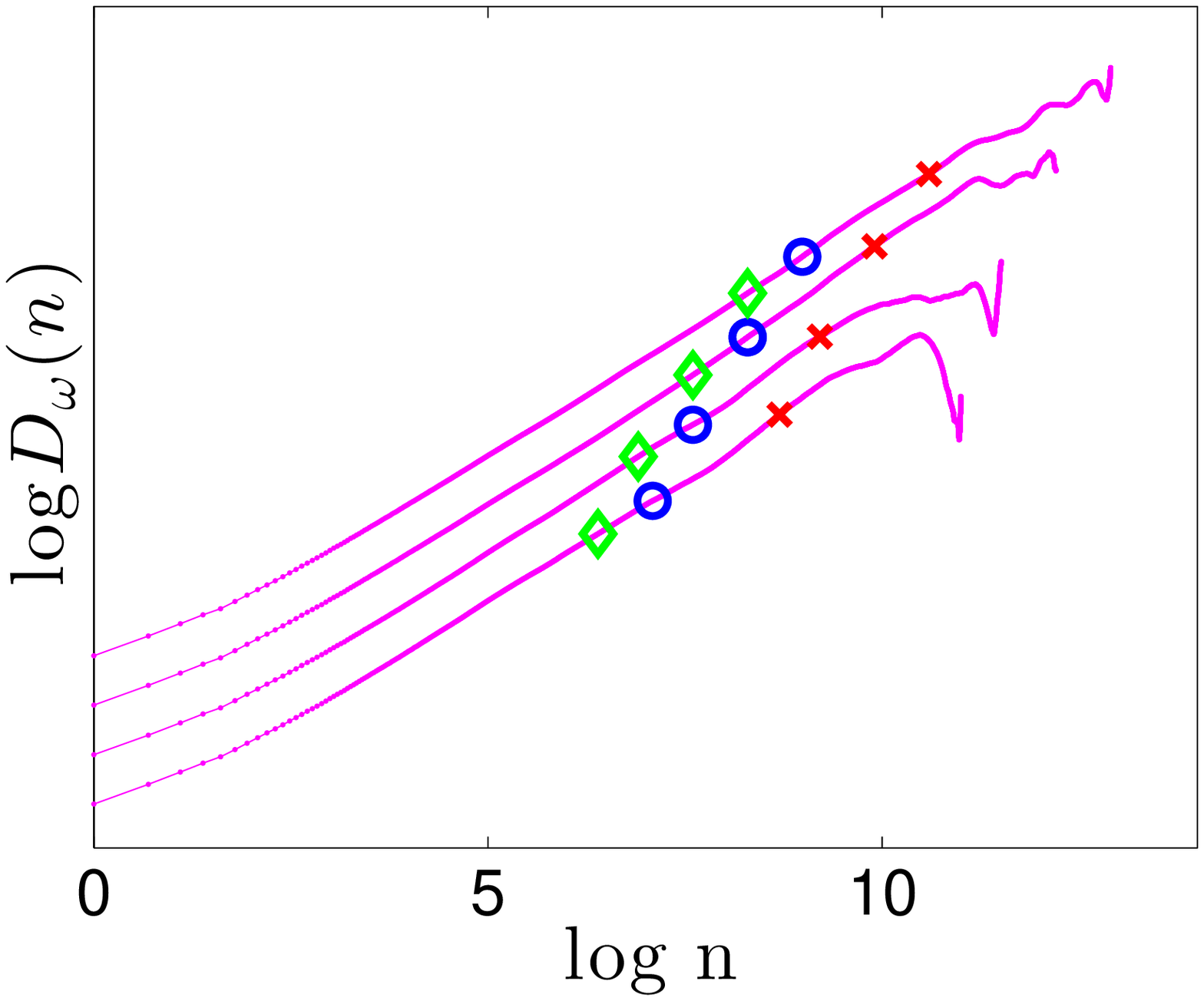}
\caption{Graph of $\log D_\omega(n)$ against $\log n$ for the Lorenz-96 system with $a=6.2$ and $\omega=0.6283$.  We take $n=1,\dots,N$ with $N=50,000$, $N=100,000$, $N=200,000$, $N=400,000$.  The various ranges used for the linear regression are marked by green diamonds for $\delta=0.01$, blue circles for $\delta=0.02$ and red crosses for $\delta=0.1$. \hl{The four graphs are spaced apart vertically so that they can be seen separately, with $N=50,000$ at the bottom, up to $N=400,000$ at the top.}}
\label{fig-L96_delta}
\end{figure}

In Figure~\ref{fig-L96_delta} we show that again a value of $\delta=0.01$ is a conservative choice for the determination of possible linear behaviour of the mean-square displacement. Note that the adequacy of $\delta=0.01$ increases with the total number of iterates $N$.

\begin{figure*}[htb]
\includegraphics[width=0.8\columnwidth]{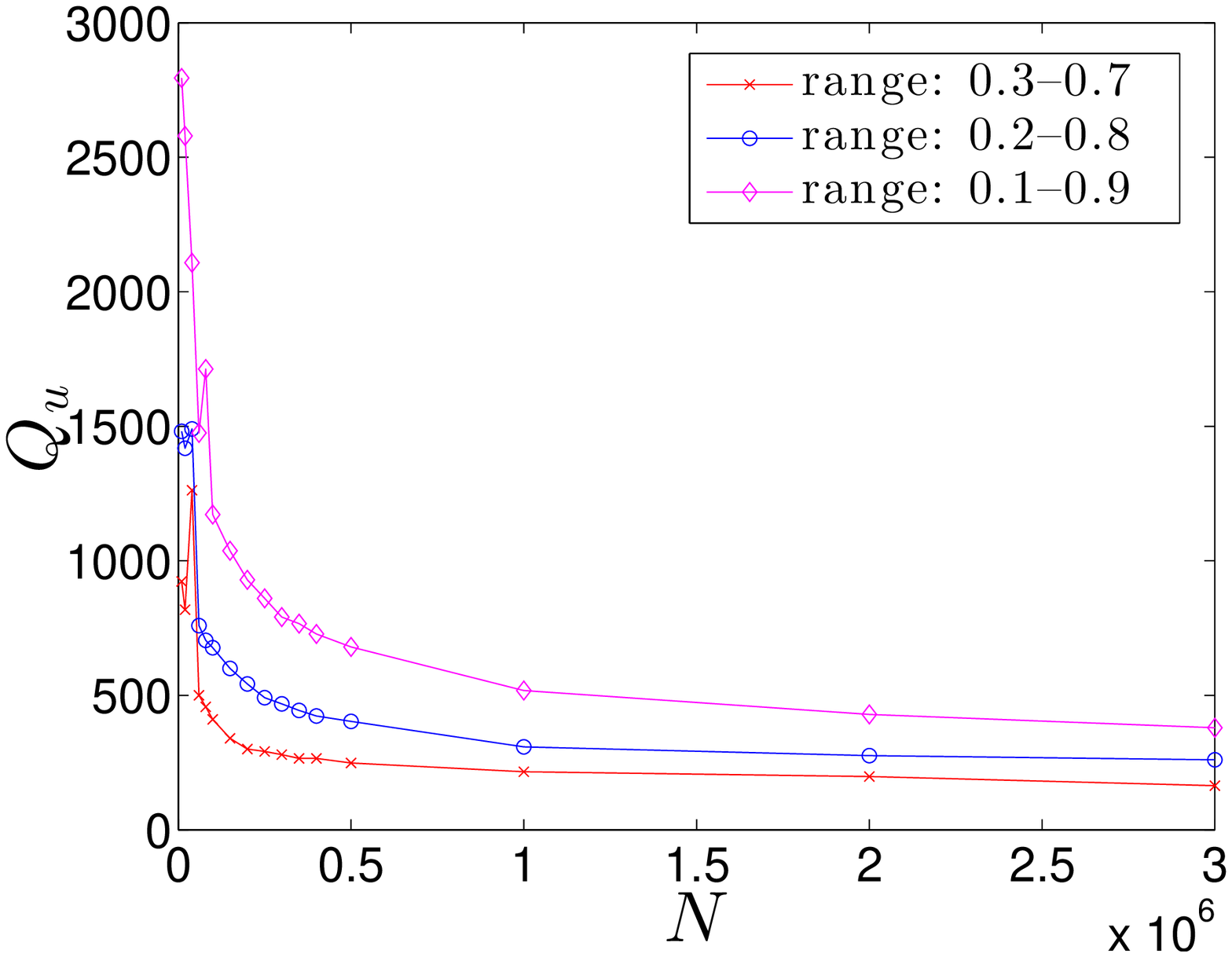}
\includegraphics[width=0.8\columnwidth]{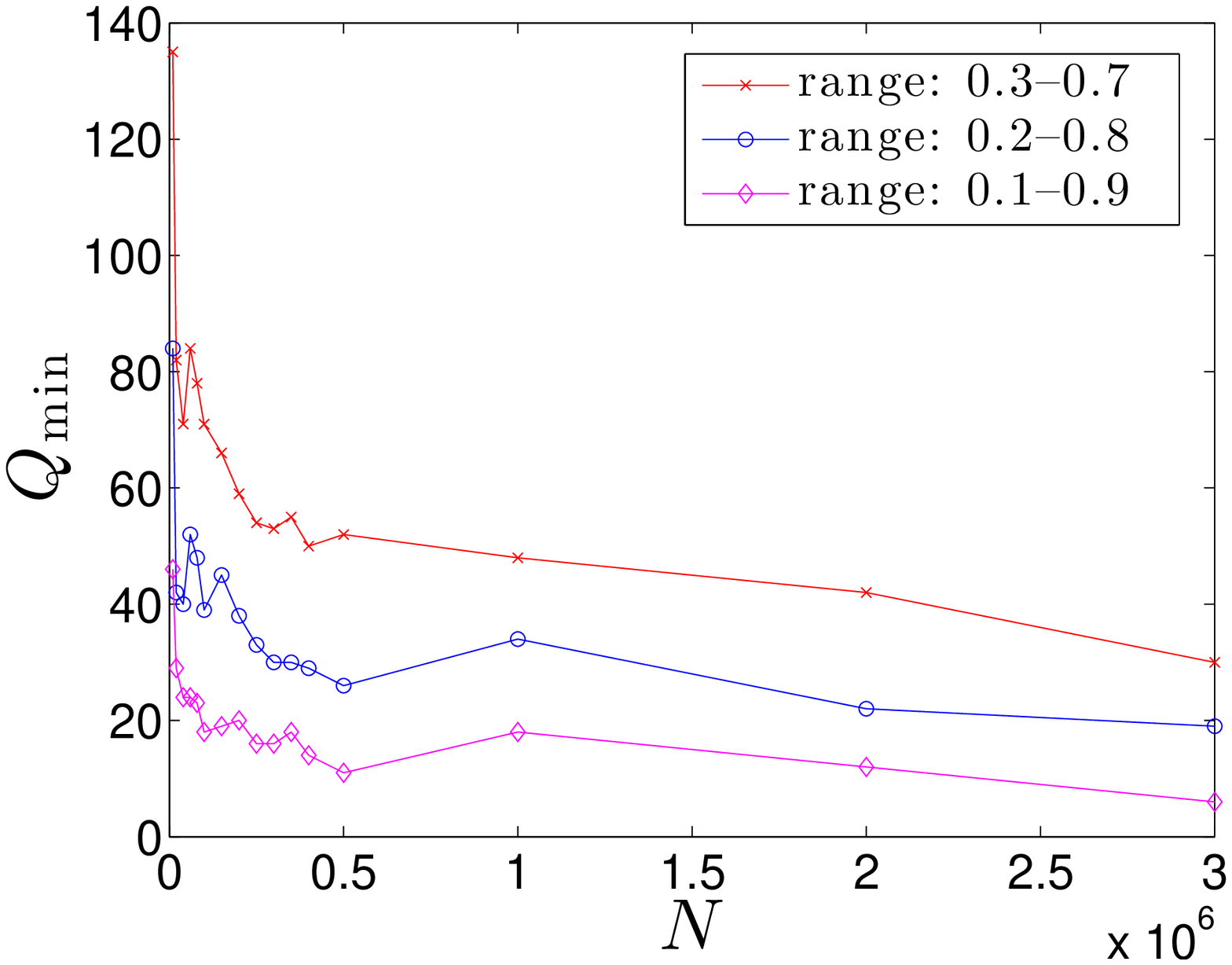}\\
\includegraphics[width=0.8\columnwidth]{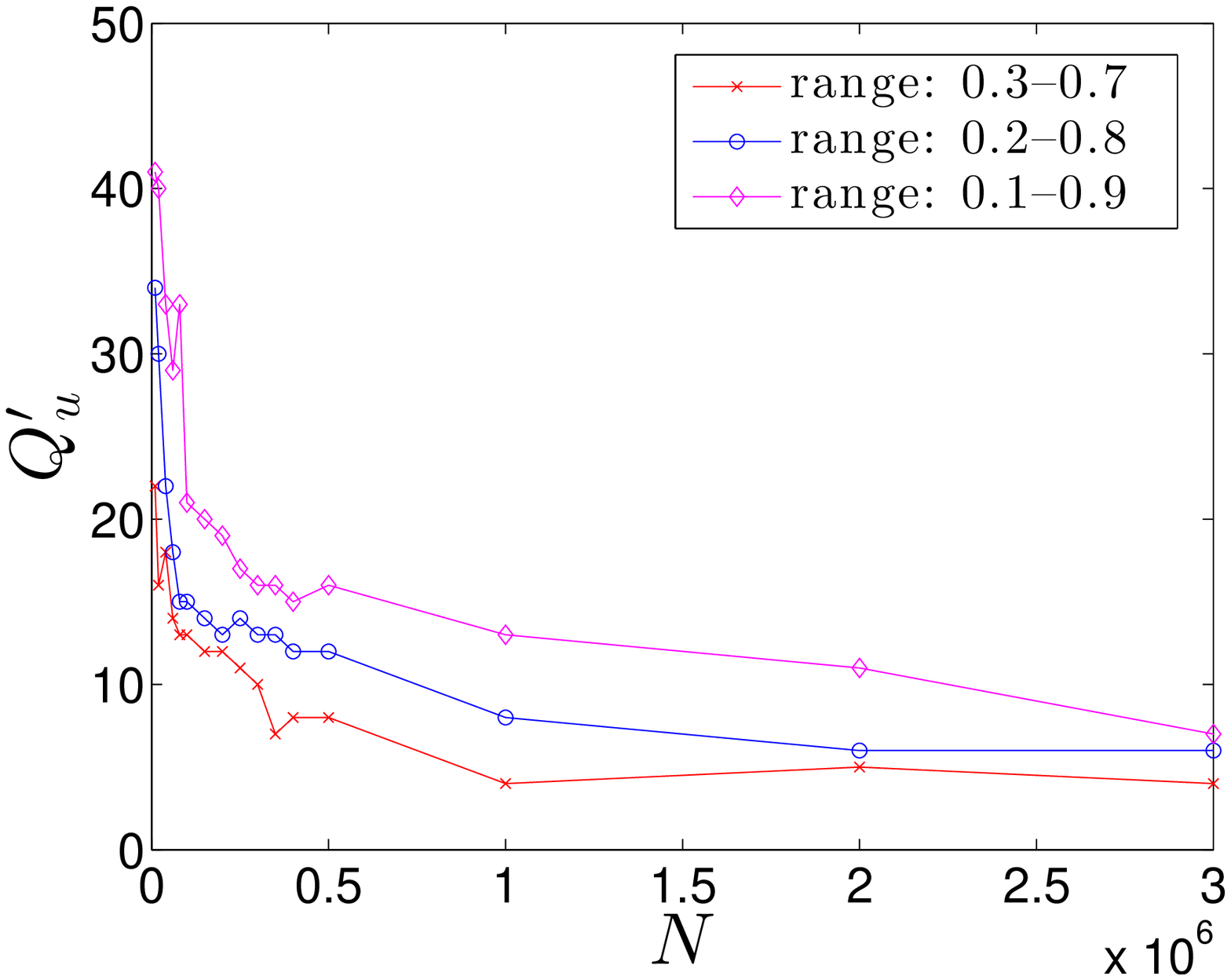}
\includegraphics[width=0.8\columnwidth]{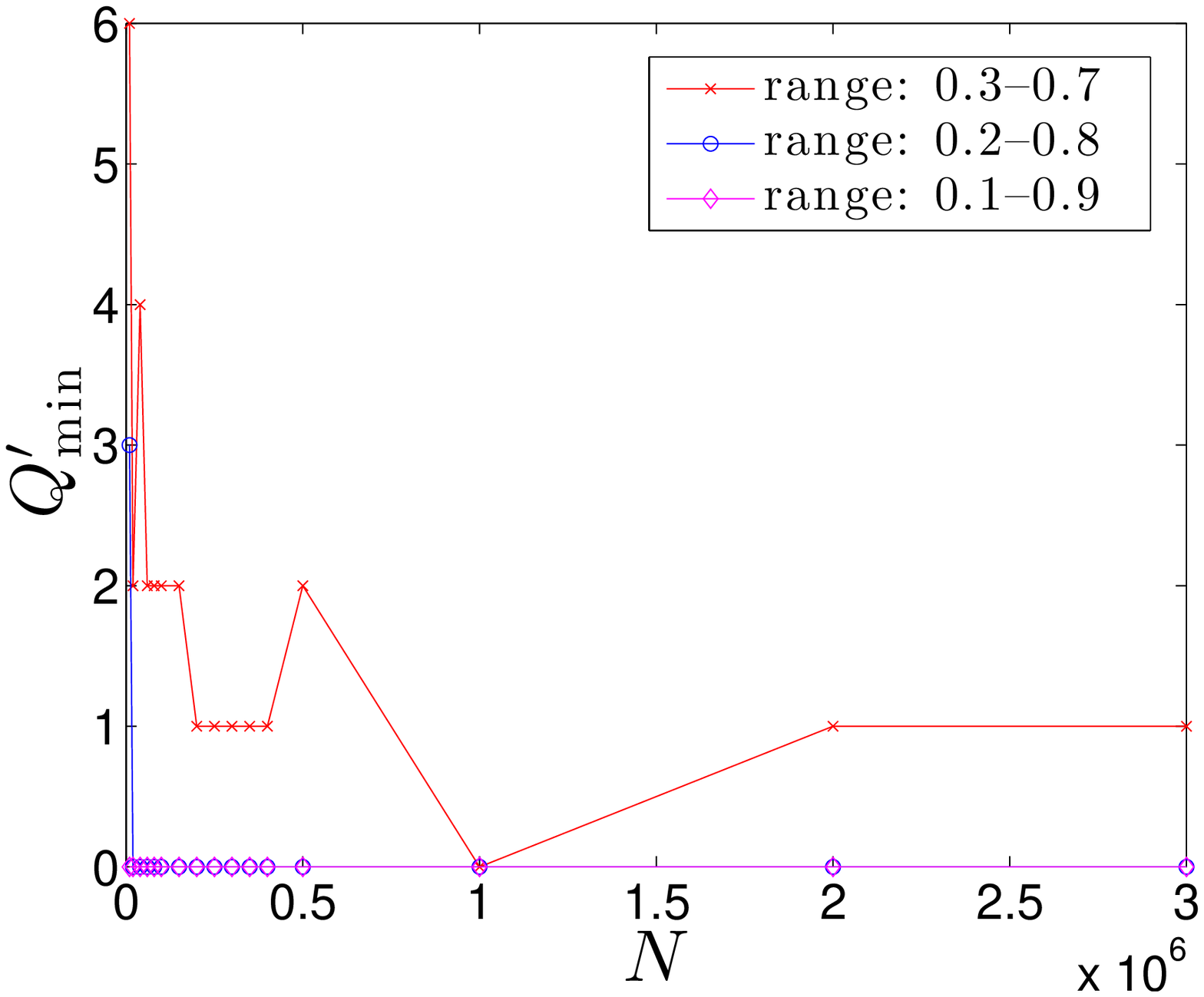}
\caption{Graphs of $Q_u$, $Q_{\rm min}$, 
$Q_u'$, $Q_{\rm min}'$ against the number of iterates $N$ for the Lorenz-96 model. 
In each case the results are shown for the three choices 
(i) $I_0=(-0.1,0.3)$, $I_1=(0.7,1.1)$, red crosses;
(ii) $I_0=(-0.1,0.2)$, $I_1=(0.8,1.1)$, blue circles;
(iii) $I_0=(-0.1,0.1)$, $I_1=(0.9,1.1)$, magenta diamonds.
}
\label{fig-L96}
\end{figure*}

Our numerical results for the Lorenz-96 model are shown in Figure~\ref{fig-L96}. The results are consistent with Conjecture~\ref{conj} that the four quantities $Q_u$, $Q_{\rm min}$, $Q'_u$, $Q'_{\rm min}$ converge to zero as $N\to\infty$. Indeed, the results in Figure~\ref{fig-L96} are comparable to those in Figure~\ref{fig-logistic} for which there was a rigorous convergence proof. The main noticeable difference is the speed of convergence since we have used $6$ times the amount of data, but that is not surprising given the increase from one dimension to $40$ dimensions and the passage from discrete to continuous time.

\begin{figure}[htb]
\includegraphics[width=0.8\columnwidth]{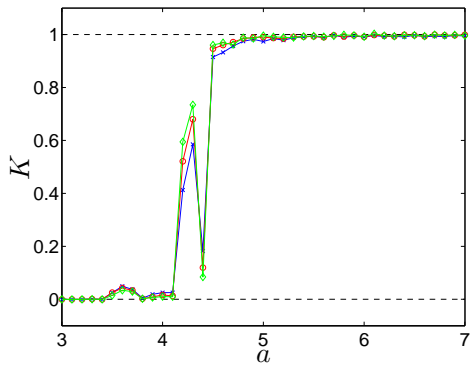}

\caption{The median value of $K(\omega)$ plotted against the parameter $a$ for the Lorenz-96 model.  
Blue crosses: $N=1,000,000$ iterates. 
Red circles: $N=2,000,000$ iterates. 
Green diamonds: $N=3,000,000$ iterates.}
\label{fig-L96_01}
\end{figure}

\begin{figure}[htb]
\includegraphics[width=0.8\columnwidth]{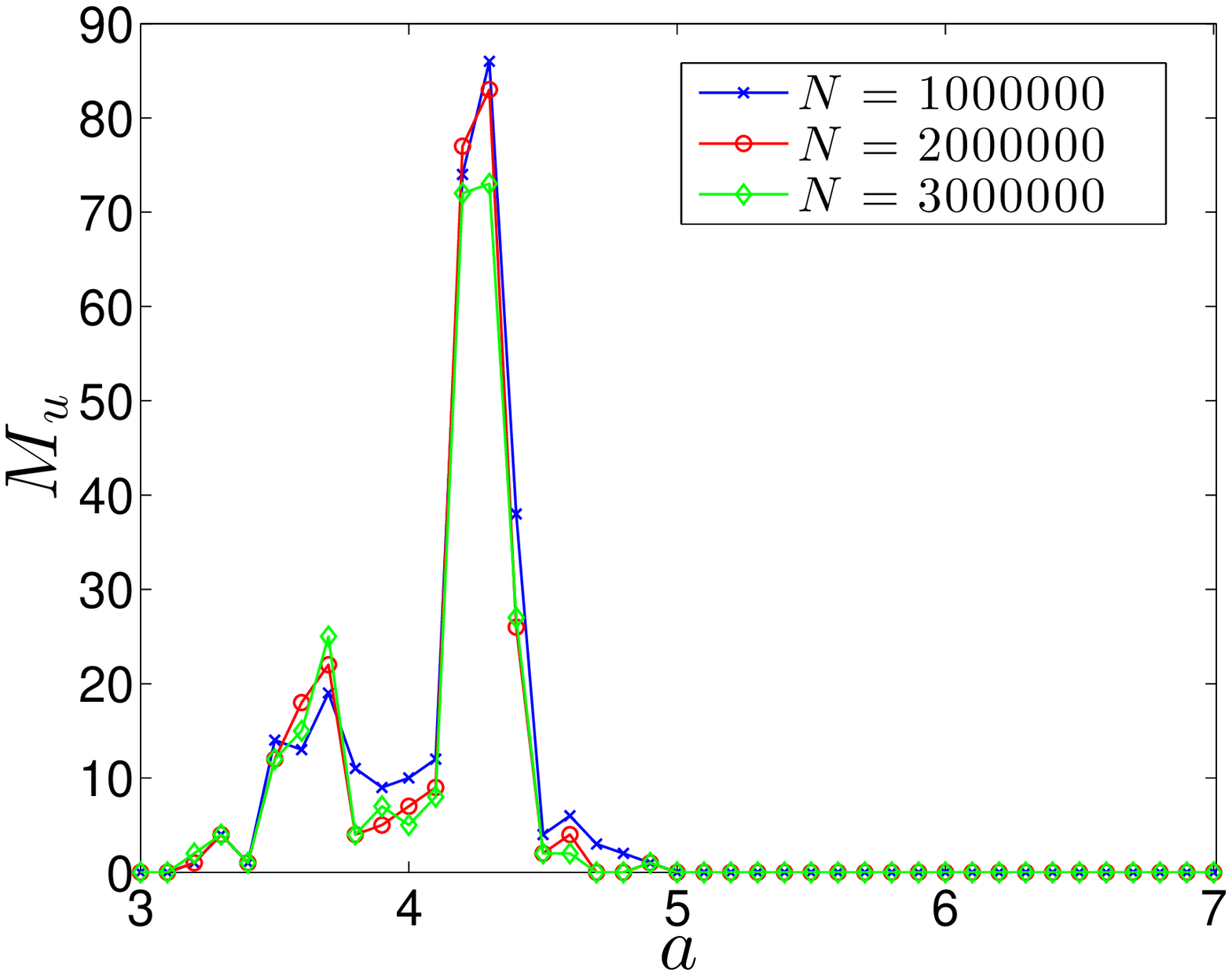}
\caption{Graph of the percentage $M_u$ of undecided values of $\omega$ plotted against the parameter $a$ 
for the Lorenz-96 model, using the range
(ii) $I_0=(-0.1,0.2)$, $I_1=(0.8,1.1)$ throughout.
Blue crosses: $N=1,000,000$ iterates. 
Red circles: $N=2,000,000$ iterates. 
Green diamonds: $N=3,000,000$ iterates.}
\label{fig-L96_outlier}
\end{figure}

Again, it is also clear that there are only a few parameter values for which the convergence is slow.  As with the logistic map, this can be compared with the corresponding results for the more quickly convergent $0$--$1$ test for chaos, see Figure~\ref{fig-L96_01}, as well as the number of outliers in Figure~\ref{fig-L96_outlier}.

\section{Summary}
\label{sec-discussion}
We have formulated a conjecture on the nature of attractors of typical dynamical systems. Our Conjecture~\ref{conj} states that typical dynamical systems are either regular, i.e.\ periodic or quasi-periodic, or strongly chaotic in the sense that they enjoy good statistical properties such as existence of SRB measures, exponential decay of correlations, large deviation principles and central limit laws.

Certain implications of Conjecture~\ref{conj} can be tested numerically and we have devised a numerical test accordingly.  The logistic map (for which a rigorous theory exists~\cite{Lyubich02}) was used as a benchmark for discussing various practical issues regarding the implementation of the test. We then proceeded to the 40-dimensional Lorenz-96 system (for which no rigorous theory is available) and showed convincing evidence that the conjecture is true also in such more complex situations.

\hl{In our numerical experiments, we have opted to fix a randomly chosen initial condition $x_1$ and then varied the parameter $a$.   
The initial condition could also be varied, thereby possibly enlarging the class of examples used for testing the conjecture.
However, for the logistic map this would not add anything new since it is known that there is a unique attractor for each value of $a$.  For the Lorenz-96 system, there is no such uniqueness result, but  since we are using the logistic map as a benchmark, it makes sense to keep the two implementations of the test as similar as possible.  (As explained in our discussion of the Palis conjecture in the introduction, it is not our aim to explore issues such as the number of coexisting attractors for a fixed parameter $a$.  Rather, we are exploring the nature of all attractors for typical initial conditions, for typical dynamical systems.)}

\hl{The main focus in this paper has been on discrete-time dissipative systems, but the conjecture applies equally to continuous time systems and to Hamiltonian systems.   For the latter, where the notion of attractor does not make sense, the conjecture would explore instead the nature of the typical asymptotic dynamics ($\omega$-limit sets) for typical initial conditions. }
 
\begin{acknowledgments}
The research of GAG was supported by the Australian Research Council; the research of IM was supported in part by the European Advanced Grant StochExtHomog (ERC AdG 320977).
\end{acknowledgments}

\bibliography{conjecture}

\end{document}